\makeatletter
  \input{jbw-switch.sty}
\makeatother

\DeclareSwitch{final}

\finaltrue

\iffinal
\documentclass[final]{itrspaper}
\else
\documentclass[draft]{itrspaper}
\fi

\usepackage{amssymb}
\usepackage{url}
\usepackage{proof}
\usepackage{color}

\newcommand{\MEM}[1]{\textup{$\mathrm{#1}$}}
\newcommand{\META}[1]{\mathit{#1}}

\def\be{\beta}

\def\l{\lambda}
\def\o{\omega}

\def\G{\Gamma}
\def\f{\rightarrow}
\def\fx{\leadsto}
\def\v{\vdash}
\def\vv{\stackrel{\circ}{\v}}
\def\<{\langle}
\def\>{\rangle}
\def\F{\displaystyle\frac}

\def\deg{\MEM{d}}
\newcommand{\mbto}{\mathbin{\to}}

\newcommand{\METArel}[0]{\META{rel}}
\newcommand{\METAfun}[0]{\META{fun}}
\newcommand{\MEMdom}[1]{\MEM{dom}(#1)}
\newcommand{\MEMran}[1]{\MEM{ran}(#1)}

\newcommand{\MEMdeg}[1]{\MEM{d}(#1)}
\newcommand{\mypair}[2]{\langle #1, #2 \rangle}
\newcommand{\func}[2]{#1 \f #2}

\iffinal
\newcommand{\mynote}[1]{}
\newcommand{\personal}[1]{}
\else
\definecolor{mygray}{gray}{0.5}
\newcommand{\mynote}[1]{
  \noindent
  \begin{footnotesize}
    \textcolor{mygray}{\textbf{[#1]}}
  \end{footnotesize}
}
\newcommand{\personal}[1]{}
\fi

\title{Realisability Semantics for Intersection Types and Expansion Variables}
\author{Fairouz Kamareddine, Karim Nour, Vincent Rahli and J. B. Wells}
\institute{http://www.macs.hw.ac.uk/ultra/}

\begin{document}

\maketitle

\begin{abstract}
  \emph{Expansion} was invented at the end of the 1970s for
  calculating \emph{principal typings} for $\lambda$-terms in type
  systems with \emph{intersection types}.  \emph{Expansion variables}
  (E-variables) were invented at the end of the 1990s to simplify and
  help mechanise expansion.  Recently, E-variables have been further
  simplified and generalised to also allow calculating type operators
  other than just intersection.  There has been much work on
  denotational semantics for type systems with intersection types, but
  none whatsoever before now on type systems with E-variables.
  Building a semantics for E-variables turns out to be challenging.
  To simplify the problem, we consider only E-variables, and not the
  corresponding operation of expansion.  We develop a realisability
  semantics where each use of an E-variable in a type corresponds to
  an independent degree at which evaluation occurs in the
  $\lambda$-term that is assigned the type.  In the $\lambda$-term
  being evaluated, the only interaction possible between portions at
  different degrees is that higher degree portions can be passed around
  but never applied to lower degree portions.  We apply this semantics
  to two intersection type systems.  We show these systems are sound,
  that completeness does not hold for the first system, and
  completeness holds for the second system when only one E-variable is
  allowed (although it can be used many times and nested).  As far as
  we know, this is the first study of a denotational semantics of
  intersection type systems with E-variables (using realisability or
  any other approach).
\end{abstract}

\section{Introduction}

Intersection types were developed in the late 1970s to type
$\lambda$-terms that are untypable with simple types; they do this by
providing a kind of finitary type polymorphism where the usage of
types is listed rather than quantified over.  They have been useful in
reasoning about the semantics of the $\lambda$-calculus, and have been
investigated for use in static program analysis.  Coppo, Dezani, and
Venneri~\cite{Cop+Dez-Cia+Ven:HBC-1980} introduced the operation of
\emph{expansion} on \emph{typings} (pairs of a type environment and a
result type) for calculating the possible typings of a term when using
intersection types.  Expansion is a crucial part of a procedure for
calculating \emph{principal typings} and thus helps support
compositional type inference.  As a simple example, the $\lambda$-term
$M=(\lambda{x}.x(\lambda{y}.yz))$ can be assigned the typing \( \Phi_1
= \langle(z : a) \vdash (((a\mbto b)\mbto b)\mbto c)\mbto c\rangle \),
which happens to be its principal typing.  The term $M$ can also be
assigned the typing \( \Phi_2 = \langle(z : a_1 \sqcap a_2) \vdash
(((a_1 \mbto b_1)\mbto b_1) \sqcap ((a_2 \mbto b_2)\mbto b_2)\mbto c
)\mbto c\rangle \), and an expansion operation can obtain $\Phi_2$
from $\Phi_1$.  Because the early definitions of expansion were
complicated, E-variables were introduced in order to make the
calculations easier to mechanise and reason about.  For example, in
System E~\cite{Car+Pol+Wel+Kfo:ESOP-2004}, the typing $\Phi_1$ from
above is replaced by \( \Phi_3 = \langle(z : e a) \vdash (e((a\mbto
b)\mbto b)\mbto c) \mbto c\rangle \), which differs from $\Phi_1$ by the
insertion of the E-variable $e$ at two places, and $\Phi_2$ can be
obtained from $\Phi_3$ by substituting for $e$ the \emph{expansion
  term} \( E = (a := a_1, b := b_1) \sqcap (a := a_2, b := b_2) \).
Carlier and Wells \cite{Car+Wel:ITRS-2004} have surveyed the history
of expansion and also E-variables.

Various kinds of denotational semantics have helped in reasoning about
the properties of entire type systems and also of specific typed
terms.  E-variables pose serious challenges for semantics.  Most
commonly, a type's semantics is given as a set of closed
$\lambda$-terms with behaviour related to the specification given by
the type.  In many kinds of semantics, the meaning of a type $T$ is
calculated by an expression $[T]_\nu$ that takes two parameters, the
type $T$ and also a valuation $\nu$ that assigns to type variables the
same kind of meanings that are assigned to types.  To extend this idea
to types with E-variables, we would need to devise some space of
possible meanings for E-variables.  Given that a type $e\,T$ can be
turned by expansion into a new type $S_1(T) \sqcap S_2(T)$, where
$S_1$ and $S_2$ are arbitrary substitutions (in fact, they can be
arbitrary further expansions), and that this can introduce an
unbounded number of new variables (both E-variables and regular type
variables), the situation is complicated.

Because it is unclear how to devise a space of meanings for expansions
and E-variables, we instead develop a space of meanings for types that
is hierarchical in the sense of having many degrees.  When assigning
meanings to types, we make each use of E-variables simply change
degrees.  We specifically avoid trying to give a semantics to the
operation of expansion, and instead treat only the E-variables.
Although this idea is not perfect, it seems to go quite far in giving
an intuition for E-variables, namely that each E-variable acts as a
kind of capsule that isolates parts of the $\lambda$-term being
analysed by the typing.  Parts of the $\lambda$-term that are typed
inside the uses of the E-variable-introduction typing rule for a
particular E-variable $e$ can interact with each other, and parts
outside $e$ can only pass the parts inside $e$ around.  The E-variable
$e$ of course also shows up in the types, and isolates the portions of
the types contributed by the portions of the term inside the
corresponding uses of E-variable-introduction.

The semantic approach we use is \emph{realisability semantics}.
Atomic types are interpreted as sets of $\lambda$-terms that are
\emph{saturated}, meaning that they are closed under $\beta$-expansion
(i.e., $\beta$-reduction in reverse).  Arrow and intersection types
are interpreted naturally by function spaces and set intersection.
Realisability allows showing \emph{soundness} in the sense that the
meaning of a type $T$ contains all closed $\lambda$-terms that can be
assigned $T$ as their result type.  This has been shown useful in
previous work for characterising the behaviour of typed
$\lambda$-terms \cite{krivine-lctm-1990}.  One also wants to show
\emph{completeness} (the converse of soundness), i.e., that every
closed $\lambda$-term in the meaning of $T$ can be assigned $T$ as its
result type.

Hindley~\cite{Hin2,Hin3,Hindley:BSTT-1997} was the first to study
completeness for a simple type system.  Then, he generalised his
completeness proof for an intersection type
system~\cite{Hindley:ISOP-1982}.  Using his completeness result for
the realisability semantics based on the sets of $\lambda$-terms
saturated by $\beta$-equivalence, Hindley has shown that simple types
are uniquely realised by the $\lambda$-terms that are typable by these
types in a type system similar to $\l_{\rightarrow}$
\cite{Barendregt:HLCS-1992} augmented with a $\beta$-equivalence rule
(this rule assigns the same typings to $\beta$-equivalent terms)
\cite{Hin2}.  He proved this result using saturation by
$\beta\eta$-equivalence w.r.t.\ a type system similar to
$\l_{\rightarrow}$ augmented with a $\beta\eta$-equivalence rule too.
Hindley also established completeness using saturation by
$\beta$-equivalence for his intersection type system
\cite{Hindley:ISOP-1982}.  In this paper, our completeness result
depends instead only on a weaker notion than $\beta$-equivalence
(saturation by $\be$-expansion).

\mynote{THE PREVIOUS TEXT WAS:
  Hindley~\cite{Hin2,Hin3,Hindley:BSTT-1997} was the first to study
  this notion of completeness for a simple type system.  Then, he
  generalised his completeness proof for an intersection type
  system~\cite{Hindley:ISOP-1982}.  Using his completeness result for
  the realisability semantics based on the sets of $\lambda$-terms
  saturated by $\beta\eta$-equivalence, Hindley has shown that simple
  types are uniquely realised by the $\lambda$-terms which are typable
  by these types.  However, this result does not hold for his
  intersection type system and he established completeness using
  saturation by $\beta\eta$-equivalence.  In this paper, our
  completeness result depends instead only on the weaker notion of
  $\beta$-equivalence. Also, we give simpler proofs that avoid needing
  $\eta$-reduction, confluence (a.k.a.\ Church-Rosser), or strong
  normalisation (SN) (we still show confluence and SN for both $\beta$
  and $\beta\eta$).}

Other work on realisability we consulted includes that by
Labib-Sami~\cite{Lab1}, Farkh and Nour~\cite{FaNo2}, and Coquand
\cite{Coquand:TLCA-2005}, although none of this work deals with
intersection types or E-variables.  Related work on realisability that
deals with intersection types includes that by Kamareddine and
Nour~\cite{kamnour}, which gives a realisability semantics with
soundness and completeness for an intersection type system.  The
system of Kamareddine and Nour is different from those in this paper,
because it allows the universal type $\omega$.  We do not know how to
build a semantics that supports both $\omega$ and E-variables.  The
method of degrees we use in this paper would need to assign $\omega$
to every degree, which is impossible.  Further work is needed on this
point.

In this paper we study the $\lambda I$-calculus typed with two
representative intersection type systems.  The restriction to $\lambda
I$ (where in $\lambda{x}.M$, the variable $x$ must be free in $M$) is
motivated by not knowing how to support the $\omega$ type.  For one of
these systems, we show that subject reduction (SR) and hence
completeness do not hold whereas for the second system, SR holds and
completeness will hold if at most one E-variable is used (although
this E-variable may be used in many places and also nested).  This is
the first paper that studies denotational semantics of intersection
type systems with E-variables, using realisability or any other
approach. One of our contributions is to outline the difficulties of
doing so.

\mynote{begin of the part added by vincent}

\personal{In this article we are not interesting in a denotational
  semantics or at least we are not interesting in an extensional
  lambda model interpreting the terms of the untyped lambda-calculus,
  but we are interesting in building a realisability semantics by
  defining sets of realisers (functions/programs satisfying the
  requirements of some specification) of types.}

The semantics we build in this paper, defines sets of realisers
(functions/programs satisfying the requirements of some specification)
of types.  Such a model can help to highlight the relation between
typable terms of the untyped lambda-calculus and types w.r.t.\ a type
system. Interpreting types in a model helps to understand the meaning
of a type (w.r.t.\ the model) which is defined as a purely syntactic
form and is clearly used as a meaningful expression\personal{a class
  of objects}. For example, the integer type, whatever its notation
is, is always used as the type of each integer. In the open problems
\mynote{Open Problems is a section in \cite{LNCS/1875}} published in
the proceedings of the Lecture Notes in Computer Science symposium
help in 1975 \cite{LNCS/1875}, it is suggested that an arrow type
expresses functionality. In that way, models based on term-models have
been built for intersection type systems
\cite{Hindley:ISOP-1982,kamnour}. In these works, intersection types
(introduced to be able to type more terms than in the Simply Typed
Lambda Calculus) are interpreted by set-theoretical intersection of
meanings. Even if expansion variables have been introduced to give a
simple formalisation of the expansion mechanism, i.e., as a syntactic
object, we are interested in the meaning of such a syntactic
object. We are particularly interested by answering these questions:
What does an expansion variable applied to a type stand for?  What are
the realisers of such a type? How can the relation between terms and
types w.r.t.\ a type system\personal{which class of objects} be
described?  How can we extend models such as the one built by
Kamareddine and Nour \cite{kamnour} to a type system with expansion?

\mynote{end of the part added by vincent}

\personal{The semantics we are considering in this article is based on
  the simple semantics \cite{Hin2} which is a particular lambda model
  (denotational model). We are interested in the study of the relation
  between typable terms and types w.r.t.\ a type system and a model.
  Moreover, interpreting types in a model helps to understand the
  meaning of a type (w.r.t.\ the model) which is defined as a purely
  syntactic form and is clearly used as a meaningful expression (as
  the integer type, whatever its notation is, which is always used as
  the type of each integer). Hence, we interpret types in the given
  model using a realisability argument (since types are interpreted by
  set of elements in the model).}

Section~\ref{Icalsec} introduces the $\l I^{\mathbb N}$-calculus,
which is the $\lambda I$-calculus with each variable marked by a
natural number \emph{degree}.  Section~\ref{sectypes} introduces the
syntax and terminology for types, and also the realisability
semantics.  Section~\ref{firstsyst} introduces our two intersection
type systems with E-variables. In one system, the syntax of types is
not restricted but in the other system it is restricted but then
extended with a subtyping relation.  We show that SR and completeness
do not hold for the first system, and that SR holds for the second
system.  We also show the soundness of the realisability semantics for
both systems and give a number of examples.  Section~\ref{complesec}
shows completeness does not hold for the second system if more than
one expansion variable is used, but does hold for a restriction of
this system to one single E-variable (which can be used in many places
and also nested).  This is an important study in the semantics of
intersection type systems with expansion variables.
Section~\ref{concsection} concludes. Full proofs can be downloaded
from the web page of the authors as well as further results that
include strong normalisation of the typable terms and the relation to
the usual unindexed $\lambda I$-calculus.

\section{The pure $\l I^{\mathbb N}$-calculus}
\label{Icalsec}

In this section we give $\l I^{\mathbb N}$, an indexed version of the
$\lambda I$-calculus where indices (which range over the set of
natural numbers ${\mathbb N} = \{0, 1, 2, \dots\}$) help categorise
the \textit{good terms} where the degree of a function is never larger
than that of its argument. This amounts to having the full $\lambda
I$-calculus at each degree (index) and creating new $\lambda I$-terms
through a mixing recipe.  Let $n, m$ be metavariables which range over
the set of natural numbers ${\mathbb N}$.  We assume that if a
metavariable $v$ ranges over a set ${\cal S}$ then $v_{i}$ for $i \geq
0$ and $v', v'',$ etc.\ also range over ${\cal S}$.  A binary relation
is a set of pairs. Let $\METArel$ range over binary relations. Let
$\MEMdom{\METArel} = \{x \mid \mypair{x}{y} \in \METArel\}$ and
$\MEMran{\METArel} = \{y \mid \mypair{x}{y} \in \METArel\}$. A
function is a binary relation $\METAfun$ such that if
$\{\mypair{x}{y}, \mypair{x}{z}\} \subseteq \METAfun$ then $y = z$.
Let $\METAfun$ range over functions. Let $\func{s}{s'} = \{\METAfun
\mid \MEMdom{\METAfun} \subseteq s \wedge \MEMran{\METAfun} \subseteq
s'\}$. We sometimes write $x : s$ instead of $x \in s$.

\begin{definition}{\ }
  \begin{enumerate}
  \item Let ${\cal V}$ be a denumerably infinite set of variables. The
    set of terms ${\cal M}$, the set of good terms ${\mathbb M}
    \subset {\cal M}$, the set of free variables $FV(M)$ of $M \in
    {\cal M}$, the degree $\deg(M)$ of a term $M$ and the joinability
    $M \diamond N$ of terms $M$ and $N$ (which ensures that in any
    term, each variable has a unique degree) are defined by
    simultaneous induction:
    \begin{itemize}
    \item If $x \in {\cal V}$, $n \in {\mathbb N}$, then $x^n \in
      {\cal M}\cap{\mathbb M}$, $FV(x^n) = \{x^n\}$, and $\deg(x^n) =
      n$.
    \item If $M, N \in {\cal M}$ such that $M \diamond N$ (see below),
      then
      \begin{itemize}
      \item $(M \, N) \in {\cal M}$, $FV((M \, N)) = FV(M) \cup FV(N)$
        and\\
        $\deg((M \; N)) = \min(\deg(M),\deg(N))$ (where $\min$ is the minimum)
      \item If $M\in {\mathbb M}$, $N\in {\mathbb M}$ and $\deg(M)
        \leq \deg(N)$ then $(M \, N)\in {\mathbb M}$.
      \end{itemize}
    \item If $M \in {\cal M}$ and $x^n \in FV(M)$, then
      \begin{itemize}
      \item $(\l x^n. M) \in {\cal M}$,
        $FV((\l x^n. M)) = FV(M)\setminus \{x^n\}$, and
        $\deg((\l x^n. M_1)) = \deg(M_1)$.
      \item If $M\in {\mathbb M}$ then $\l x^n. M\in {\mathbb M}$.
      \end{itemize}
    \end{itemize}

  \item Let $M, N \in {\cal M}$.  We say that $M$ and $N$ are joinable
    and write $M \diamond N$ iff $\forall x\in {\cal V}$, if $x^m \in
    FV(M)$ and $x^n \in FV(N)$, then $m = n$.  If ${\cal X} \subseteq
    {\cal M}$ such that $\forall M, N \in {\cal X}, M \diamond N$, we
    write, $\diamond{\cal X}$.  If ${\cal X} \subseteq {\cal M}$ and
    $M \in {\cal M}$ such that $\forall N \in {\cal X}, M \diamond N$,
    we write, $ M \diamond{\cal X}$.

  \item We adopt the usual definition
    \cite{Barendregt:LCSS-1984,krivine-lctm-1990} of subterms and the
    convention for parentheses and their omission.  Note that a
    subterm of $M \in {\cal M}$ (resp.\ ${\mathbb M}$) is also in
    ${\cal M}$ (resp.\ ${\mathbb M}$).  We let $x, y, z,$ etc.\ range
    over ${\cal V}$ and $M, N, P,$ etc.\ range over ${\cal M}$ and use
    $=$ for syntactic equality.

  \item For each $n \in {\mathbb N}$, we let: \hspace{0.2in} $\bullet$
    ${\cal M}^n = \{M \in {\cal M} \mid \deg(M) =  n\}$\\
    $\bullet$ ${\cal M}^{> n} = {\cal M}^{\geq n+1}$ \hspace{0.02in}
    $\bullet$ ${\cal M}^{\geq n} = \{M \in {\cal M} \mid \deg(M) \geq
    n\}$ \hspace{0.02in} $\bullet$ ${\mathbb M}^n = {\mathbb
      M}\cap{\cal M}^n$

  \item  For $m \geq 0$, $M[(x^{n_i}_i:=N_i)_{1 \leq i
      \leq m}]$ (or simply $M[(x^{n_i}_i:=N_i)_m]$), the simultaneous
    substitution of $N_i$ for all free occurrences of $x^{n_i}_i$ in
    $M$ only matters when $\diamond {\cal X}$ where ${\cal X} =
    \{M\}\cup\{N_i \mid 1 \leq i \leq m\} \subseteq {\cal M}$. Hence we
    restrict substitution accordingly to incorporate the $\diamond$
    condition. With ${\cal X}$ as above, $M[(x^{n_i}_i:=N_i)_m]$ is
    only defined when $\diamond {\cal X}$. We write
    $M[(x^{n_i}_i:=N_i)_{1 \leq i \leq 1}]$ as $M[x^{n_1}_1:=N_1]$.

  \item We take terms modulo \textit{$\alpha$-conversion} given by:
    $\l x^n. M = \l y^n. (M[x^n:=y^n]) \mbox{ where } \forall m, y^m
    \not \in FV(M)$.  We use the Barendregt convention (BC) where the
    names of bound variables differ from the free ones and where we
    rewrite terms so that not both $\l x^n$ and $\l x^m$ co-occur when
    $n\not =m$.

  \item A relation $R$ on ${\cal M}$ is \textit{compatible} iff for
    all $M, N, P \in {\cal M}$:
    \begin{itemize}
    \item If $\mypair{M}{N} \in R$ and $x^n \in FV(M)\cap FV(N)$ then
      $\mypair{\l x^n. M}{\l x^n. N} \in R$.
    \item If $\mypair{M}{N} \in R$, $M \diamond P$ and $N \diamond P$ then
      $\mypair{MP}{NP} \in R$ and $\mypair{PM}{PN} \in R$.
    \end{itemize}

  \item The reduction relation $\rhd_\be $ on ${\cal M}$ is defined as
    the least compatible relation closed under the rule: $(\l x^n. M)N
    \rhd_\be M[x^n:=N]$ if $\deg(N) = n.$

  \item We denote by $\rhd_\beta^*$ the reflexive and transitive
    closure of $\rhd_\beta$.  We denote by $\simeq_{\beta}$ the
    equivalence relation induced by $\rhd_{\beta}^*$.
  \end{enumerate}
\end{definition}

Beta reduction is well defined on the $\l I^{\mathbb N}$-calculus,
i.e., if $M \in {\cal M}$ and $M\rhd_\beta N$ then $N \in {\cal M}$.
(Note that because $\MEMdeg{x^{0}} = 0 \not = 1 = \MEMdeg{z^{1}}$,
then $(\l x^{0}. x^{0}y^{0})z^{1} \not {\rhd_{\be}} z^{1}y^{0}$.)
Hence, $\rhd_\beta^*$ is also well defined on ${\cal M}$.  Beta
reduction preserves the free variables, degrees and goodness of terms,
i.e., if $M \rhd_\beta^* N$ then $FV(M) = FV(N)$, $\deg(M) = \deg(N)$
and $M$ is good iff $N$ is good.

The next definition turns terms of degree $n$ into terms of higher
degrees and also, if $n >0$, they can be turned into terms of lower
degrees. Note that $^+$ and $^-$ are well behaved operations with
respect to all that matters (free variables, reduction, joinability,
substitution, etc.).

\begin{definition}{\ }
  \label{def:plusandminus}
  \begin{enumerate}
  \item We define $^+: {\cal M} \f {\cal M}$ and $^-: {\cal M}^{>
      0} \f {\cal M}$ by:\\
    $\bullet$ $(x^n)^+ = x^{n+1}$ \hspace{0.07in} $\bullet$ $(M_1 \;
    M_2)^+ = M_1^+ \;
    M_2^+$ \hspace{0.07in} $\bullet$ $(\l x^n.M)^+ = \l x^{n+1}.M^+$\\
    $\bullet$ $(x^n)^- = x^{n-1}$ \hspace{0.07in} $\bullet$ $(M_1 \;
    M_2)^- = M_1^- \; M_2^-$ \hspace{0.07in} $\bullet$ $(\l x^n.M)^- =
    \l x^{n-1}.M^-$
  \item Let ${\cal X} \subseteq {\cal M}$.  If $\forall M \in {\cal
      X}$, $\deg(M)>0$, we write $\deg({\cal X}) >0$.  We define:\\
    $\bullet$ ${\cal X}^+ = \{M^+ \mid M \in {\cal X}\}$
    \hspace{0.4in} $\bullet$ If $\deg({\cal X}) >0$, ${\cal X}^- =
    \{M^- \mid M \in {\cal X}\}$.
  \item We define $M^{-n}$ by induction on $\deg(M) \geq n \geq 0$. If
    $n = 0$ then $M^{-n} = M$ and if $n \geq 0$ then $M^{-(n+1)} =
    (M^{-n})^-$.
  \end{enumerate}
\end{definition}

\section{The types and their realisability semantics}
\label{sectypes}

This paper studies two type systems.  In the first, there are
no restrictions on where the arrow occurs. In the second, arrows
cannot occur to the left of intersections or expansions.  The next
definition gives these two basic sets of types and the notions of
a degree of a type and of a good type.

\begin{definition}[Types, good types, degree of a type]
  \label{degdef}\mbox{} \hfill
  \begin{enumerate}
  \item Assume two denumerably infinite sets ${\cal A}$ (atomic types)
    and ${\cal E}$ (expansion variables).  Let $a, b, c,$ etc.\ range
    over ${\cal A}$ and $e$ range over ${\cal E}$.
  \item The sets of types ${\cal T}$, ${\mathbb U}$ and ${\mathbb T}$
    are defined by ${\cal T}::= {\cal A} \mid {\cal T} \f {\cal T}
    \mid {\cal T}
    \sqcap {\cal T}\mid {\cal E}{\cal T}$ and\\
    ${\mathbb U} ::= {\mathbb U} \sqcap {\mathbb U} \mid {\cal
      E}{\mathbb U} \mid {\mathbb T}$ \hspace{0.01in} where
    \hspace{0.01in} ${\mathbb T} ::= \; {\cal A} \mid {\mathbb U} \f
    {\mathbb T}$ (note that ${\mathbb T}$ and ${\mathbb U}$ are
    defined simultaneously).  Note that ${\mathbb T} \subseteq
    {\mathbb U} \subseteq {\cal T}$.  We let $T, U, V, W$ (resp.\ $T$,
    resp.\ $U, V, W$) range over ${\cal T}$ (resp.\ ${\mathbb T}$,
    resp.\ ${\mathbb U}$).  We quotient types by taking $\sqcap$ to be
    commutative,
    associative,
    idempotent,
    and to satisfy $e(U_1 \sqcap U_2) = eU_1 \sqcap eU_2$.

  \item Denote $e_{i_l} \dots e_{i_n}$ by $\vec{e}_{i(l:n)}$ and
    $U_n\sqcap U_{n+1} \dots \sqcap U_{m}$ by $\sqcap_{i=n}^{m}U_i$
    ($n \leq m$).

  \item We define a function $\deg: {\cal T} \f {\mathbb N}$ by
    (hence $\deg$ is also defined on ${\mathbb U}$):\\
    $\bullet \; \deg(a) = 0$\hspace{1in}
    $\bullet \; \deg(U \f T) = \min(\deg(U),\deg(T))$\\
    $\bullet \; \deg(eU) = \deg(U)+1$\hspace{0.39in} $\bullet \;
    \deg(U \sqcap V) = \min(\deg(U),\deg(V))$.

  \item We define the good types on ${\cal T}$ by (this also defines
    good types on ${\mathbb U}$):\\
    $\bullet$ If $a \in {\cal A}$, then $a$ is good \hspace{0.07in}
    $\bullet$ If $U$ is good and $e \in {\cal E}$, then $eU$ is good\\
    $\bullet$ If $U, T$ are good and $\deg(U) \geq \deg(T)$, then $U
    \f T$ is good\\
    $\bullet$ If $U, V$ are good and $\deg(U) = \deg(V)$, then $U
    \sqcap V$ is good
  \end{enumerate}
\end{definition}

\begin{definition}[Environments]
  \label{envdef}
  \begin{enumerate}
  \item A type environment is a set $\{x^{n_i}_i:U_i \mid 1 \leq i \leq n
    \mbox{ where } n \geq 0 \mbox{ and } \forall 1 \leq i, j \leq n,
    \mbox{ if } i \not = j \mbox { then } x^{n_i}_i \not =
    x^{n_j}_j\}$. We denote such environment (call it $\G$) by
    $x^{n_1}_1:U_1, x^{n_2}_2:U_2, \dots , x^{n_n}_n:U_n$ or simply by
    $(x^{n_i}_i:U_i)_n$ and define $dom(\G) = \{x^{n_i}_i \mid 1 \leq i
    \leq n\}$.  We use $\G, \Delta$ to range over
    environments and write $()$ for the empty environment.\\
    Of course on ${\cal T}$, type environments take variables in
    ${\cal V}$ to ${\cal T}$.  On $ {\mathbb U}$, they take variables
    in ${\cal V}$ to $ {\mathbb U}$.
    \begin{itemize}
    \item We say that $\G$ is good iff , for every $1 \leq i \leq k$, $U_i$ is
      good.
    \item We say that $\deg(\G) > 0$ iff for every $1 \leq i \leq k$, $\deg(U_i)>
      0$ and $n_i > 0$.
    \end{itemize}

  \item If $\G = (x^{n_i}_i:U_i)_n$ and $x^m \not \in dom(\G)$, then
    we write $\G, x^m : U$ for the type environment $x^{n_1}_1:U_1,
    \dots , x^{n_n}_n:U_n , x^m : U$.
  \item Let $\G_1 = (x^{n_i}_i:U_i)_n,(y^{m_j}_j:V_j)_m$ and $\G_2 =
    (x^{n_i}_i:U'_i)_n,(z^{r_k}_k:W_k)_r$. We write $\G_1 \sqcap \G_2$
    for the type environment $(x^{n_i}_i:U_i \sqcap
    U'_i)_n,(y^{m_j}_j:V_j)_m,(z^{r_k}_k:W_k)_r$. Note that $dom(\G_1
    \sqcap \G_2) = dom(\G_1) \cup dom(\G_2)$ and that $\sqcap$ is
    commutative, associative and idempotent on environments.
  \item $e\G = (x^{n_i+1}_i:eT_i)_n$ where $\G = (x^{n_i}_i:T_i)_n$.
    So $e(\G_1\sqcap\G_2) = e\G_1\sqcap e\G_2$.
  \item We say that $\G_1$ is joinable with $\G_2$ and write
    $\G_1 \diamond \G_2$ iff\\
    $\phantom\quad\quad\phantom \quad \forall x\in {\cal V}$, if $x^m
    \in dom(\G_1)$ and $x^n \in dom(\G_2)$, then $m = n$.
  \end{enumerate}
\end{definition}

\begin{definition}[Degree decreasing of a type]
  \begin{enumerate}
  \item If $\deg(U) > 0$, we inductively define the type $U^-$ by:
    $\bullet$ $(U_1 \sqcap U_2)^- = U_1^- \sqcap U_2^-$ \hspace{0.5in}
    \hspace{0.4in} $\bullet$$(eU)^- = U$\\
    If $\deg(U) \geq n \geq 0$, $U^{-n}$ is defined as for $M^{-n}$ in
    definition~\ref{def:plusandminus}.
  \item If $\G = (x^{n_i}_i : U_i)_k$ and $\deg(\G) > 0$, then
    we let $\G^- = (x^{n_i-1}_i : U^-_i)_k$.\\
    If $\deg(\G) \geq n \geq 0$, $\G^{-n}$ is defined as for $M^{-n}$ in
    definition~\ref{def:plusandminus}.
  \item If $U$ is a type and $\G$ is a type environment such that
    $\deg(\G) > 0$ and $\deg(U)> 0$, then we let $(\<\G \v_2 U\>)^- =
    (\<\G^- \v_2 U^-\>)$.
  \end{enumerate}
\end{definition}

Saturated sets
and the interpretations and meanings of types are crucial to a realisability semantics:

\begin{definition}[Saturated sets]
  \label{setdefs}
  Let ${\cal X},{\cal Y} \subseteq {\cal M}$.
  \begin{enumerate}
  \item We use ${\cal P}({\cal X})$ to denote the powerset of ${\cal
      X}$, i.e.\ $\{{\cal Y} \mid {\cal Y} \subseteq {\cal X}\}$.
  \item We let ${\cal X} \fx {\cal Y} = \{M \in {\cal M} \mid \forall N
    \in {\cal X}$, if $M \diamond N$ then $M \; N \in {\cal Y}\}$.
  \item ${\cal X}$ is saturated iff whenever $M \rhd_\be^*
    N$ and $N \in {\cal X}$, then $M \in {\cal X}$.
  \end{enumerate}
\end{definition}

\begin{definition}[Interpretations and meaning of types]
  \label{intdef}
  Let ${\cal V} = {\cal V}_1\cup {\cal V}_2$ where ${\cal V}_1\cap
  {\cal V}_2 = \emptyset$ and ${\cal V}_1, {\cal V}_2$ are both
  denumerably infinite.
  \begin{enumerate}
  \item Let $x \in {\cal V}_1$ and $n \in {\mathbb N}$.  We define
    ${\cal N}_x^n = \{x^n \; N_1...N_k \in {\mathbb M} \mid k \geq
    0\}$.
  \item An interpretation ${\cal I}: {\cal A} \f {\cal P}({\cal
      M}^0)$ is a function such that for all $a \in {\cal A}$:\\
    $\bullet \; {\cal I}(a)$ is saturated \hspace{0.4in} and
    \hspace{0.4in} $\bullet \; \forall x \in {\cal V}_1, \; {\cal
      N}_x^0 \subseteq {\cal I}(a)\subseteq {\mathbb M}^0$.

  \item Let an interpretation ${\cal I}: {\cal A} \f {\cal
      P}({\cal M}^0)$. We extend ${\cal I}$ to ${\cal
      T}$ (hence this includes ${\mathbb U}$)
    as follows: $\; \; \bullet \; {\cal I}(eU) = {\cal I}(U)^+$ \hspace{0.5in}
    $\bullet \; {\cal I}(U \sqcap V) = {\cal I}(U) \cap {\cal I}(V)$
    \hspace{0.5in} $\bullet \; {\cal I}(U \f T) = {\cal I}(U) \fx
    {\cal I}(T)$

    Because $\cap$ is commutative, associative, idempotent, and $({\cal
    X} \cap {\cal Y})^+ = {\cal X}^+ \cap {\cal Y}^+$, ${\cal I}$
    is well defined.
  \item Let $U \in {\cal T}$ (hence $U$ can be in ${\mathbb U}$).  We
    define the meaning $[U]$ of $U$ by:\\
      $[U]= \{M \in {\cal M} \mid M \mbox{ is closed and } M \in
      \bigcap_{{\cal I}\mbox{ interpretation }} {\cal I}(U)\}$.
  \end{enumerate}
\end{definition}

It is easy to show that if $x^n \; N_1...N_k \in {\cal N}_x^n$ then
$\forall~1 \leq i \leq k$, $\deg(N_i) \geq n$.

Type interpretations are saturated and interpretations of good types
contain only good terms.

\section{The typing systems $\v_1$ and $\v_2$}
\label{firstsyst}

In this section we introduce $\v_1$ and $\v_2$, our two intersection
type systems with expansion variables.  In $\v_1$, types are not
restricted and SR fails.  In $\v_2$, the syntax of types is restricted
in the sense that arrows cannot occur to the left of intersections or
expansions. In order to guarantee SR for this type system (and hence
completeness later on), we introduce a subtyping relation which will
allow intersection type elimination (something not available in the
first type system).

\begin{definition}
  Let $i \in \{1,2\}$.
  The type system $\v_1$ (resp.\ $\v_2$) uses the set ${\cal T}$
  (resp.\ ${\mathbb U}$) of definition~\ref{degdef}.  We follow
  \cite{Car+Wel:ITRS-2004} and write type judgements as $M: \<\G \v
  U\>$ instead of the traditional format of $\G \v M: U$.  The typing
  rules of $\v_i$ are (recall that when used for $\v_1$, $U$ and $T$
  range over ${\cal T}$, and when used for $\v_2$, $U$ ranges over
  ${\mathbb U}$ and $T$ ranges over ${\mathbb T}$) of
  figure~\ref{fig:typesystemrules} (left).  In the last clause, the
  binary relation $\sqsubseteq$ is defined on ${\mathbb U}$ by the
  rules of figure~\ref{fig:typesystemrules} (right).

  Let $\Phi$ denote types in ${\mathbb U}$, or
  environments $\G$ or typings $\<\G\v_2U\>$. When $\Phi
  \sqsubseteq \Phi'$, then $\Phi$ and $\Phi'$ belong to the same set
  (${\mathbb U}$/environments/typings).
  Let $\G$ be an environment, $U \in {\cal T}$ and $M \in {\cal
    M}$.\\
  $\bullet$ We say that $\G$ is $\v_i$-legal iff there are $M, U$ such
  that $M : \< \G \v_i U\>$.\\
  $\bullet$ We say that $\<\G \v_i U\>$ is good iff $\G$ and $U$
  are good.\\
  $\bullet$  We say that $\deg(\<\G \v_i U\>) > 0$ iff $\deg(\G) > 0$ and
  $\deg(U)> 0$.
\end{definition}
{\footnotesize
\begin{figure}[t]
  \begin{tabular}{|c|c|}
    \hline
    \begin{small}
      \begin{tabular}{c}
        \infer[(ax)]{x^n : \<(x^n:T) \v_1 T\>}{T \mbox{ good} & \deg(T)=
          n}\\
        \\
        \infer[(ax)]{x^0 : \<(x^0:T) \v_2 T\>}{T \mbox{ good}}\\
        \\
        \infer[(\f_I)]{\l x^n. M : \<\G \v_i U \f T\>}{M : \<\G,(x^n:U)
          \v_i T\>}\\
        \\
        \infer[(\f_E)]{M_1 M_2 : \<\G_1 \sqcap \G_2 \v_i T\>}{M_1 :
          \<\G_1 \v_i U \f T\> & M_2 : \<\G_2 \v_i U\> & \G_1 \diamond
          \G_2}\\
        \\
        \infer[(\sqcap)]{M : \<\G_1 \sqcap \G_2 \v_i U_1 \sqcap U_2\>}{M:
          \<\G_1 \v_i U_1\> & M : \<\G_2 \v_i U_2\>}\\
        \\
        \infer[(exp)]{M^+ : \< e\G \v_i eU\>}{M : \<\G \v_i U\>}\\
        \\
        \infer[(\sqsubseteq)]{M : \<\G' \v_2 U'\>}{M : \<\G \v_2 U\> &
          \<\G \v_2 U\> \sqsubseteq \<\G' \v_2 U'\>}
      \end{tabular}
    \end{small}
    &
    \begin{small}
      \begin{tabular}{c}
        \infer[(ref)]{\Phi \sqsubseteq \Phi}{}\\
        \\
        \infer[(tr)]{\Phi_1 \sqsubseteq \Phi_3}{\Phi_1 \sqsubseteq \Phi_2 &
          \Phi_2 \sqsubseteq \Phi_3}\\
        \\
        \infer[(\sqcap_e)]{U_1 \sqcap U_2 \sqsubseteq U_1}{U_2 \mbox{ good} &
          \deg(U_1) = \deg(U_2)}\\
        \\
        \infer[(\sqcap)]{U_1 \sqcap U_2 \sqsubseteq V_1 \sqcap V_2}{U_1
          \sqsubseteq V_1 & U_2 \sqsubseteq V_2}\\
        \\
        \infer[(\f)]{U_1 \f T_1 \sqsubseteq U_2 \f T_2}{U_2 \sqsubseteq U_1
          & T_1 \sqsubseteq T_2}\\
        \\
        \infer[(\sqsubseteq_{exp})]{eU_1 \sqsubseteq eU_2}{U_1 \sqsubseteq
          U_2}\\
        \\
        \infer[(\sqsubseteq_c)]{\G, (y^n : U_1) \sqsubseteq \G, (y^n :
          U_2)}{U_1 \sqsubseteq U_2}\\
        \\
        \infer[(\sqsubseteq_{\<\>})]{\<\G_1 \v_2 U_1\> \sqsubseteq \<\G_2
          \v_2 U_2\>}{U_1 \sqsubseteq U_2 & \G_2 \sqsubseteq \G_1}
      \end{tabular}
    \end{small}\\
    \hline
  \end{tabular}
  \label{fig:typesystemrules}
  \vspace{-0.1in}
  \caption{Typing rules / Subtyping rules}
\end{figure}
}
We show that typable terms are good, have good types, and have the
same degree as their types and that all legal contexts are good.  We
also show that no $\beta$-redexes are blocked in a typable term.



SR for $\beta$ using $\v_1$ fails:
let $a, b, c$ be different elements of ${\cal A}$.  Although
$(\l x^0.  x^0 x^0) (y^0 z^0) \rhd_{\beta} (y^0 z^0)(y^0 z^0)$ and
$(\l x^0. x^0 x^0) (y^0 z^0) : \< y^0 : b \f ((a\f c)\sqcap a), z^0: b
\v_1 c\>$, it is not possible that $(y^0 z^0)(y^0 z^0) : \< y^0 : b \f
((a\f c)\sqcap a), z^0: b \v_1 c\>$.

Nevertheless, we show that SR and subject expansion for
$\beta$ using $\v_2$ holds.  This will be used in the proof of completeness (more specifically
in lemma~\ref{comp} which is basic for the completeness
theorem~\ref{comple}).

\begin{lemma}[Subject reduction and expansion for $\beta$]
  \label{betaetasr}
  \label{subject-expansion-star}
\mbox{} \hfill
  \begin{enumerate}
  \item If $M : \<\G \v_2 U\>$ and $M \rhd^*_{\beta} N$, then $N :
    \<\G \v_2 U\>$.
  \item If $N : \<\G \v_2 U\>$ and $M \rhd_\be^* N$ then $M : \<\G
    \v_2 U\>$.
  \end{enumerate}
\end{lemma}

The semantics given in section~\ref{sectypes} is sound with respect to
$\v_1$ and $\v_2$, because if ${\cal I}$ is an interpretation and $U
\sqsubseteq V$ then ${\cal I}(U) \subseteq {\cal I}(V)$.

\begin{lemma}[Soundness of $\v_1$/$\v_2$]
  \label{adeq}
  \label{adeq'}
  Let $i \in \{1,2\}$, ${\cal I}$ be an interpretation, $M :
  \<(x^{n_j}_j:U_j)_n \v_i U\>$ and $\forall 1 \leq j \leq n$, $N_j
  \in {\cal I}(U_j)$.  If $M[(x^{n_j}_j:=N_j)_n] \in {\cal M}$, then
  $M[(x^{n_j}_j:=N_j)_n] \in {\cal I}(U)$.
\end{lemma}

Hence, if $M : \<() \v_i U\>$, then $M \in [U]$.  The next lemma puts
the realisability semantics in use.

\begin{lemma}
  \label{exlem}
  \label{newexlem}
  \begin{enumerate}
  \item \label{exseven} \label{newexseven} $[(a \sqcap b) \f a] = \{M
    \in {\mathbb M}^0 \mid M \rhd_\be^*\l y^0.y^0\}$.
  \item \label{exeight} It is not possible that $\l y^0.y^0 :\<()\v_1
    (a \sqcap b) \f a\>$.
  \item \label{newexeight} $\l y^0.y^0 :\<()\v_2 (a \sqcap b) \f a\>$.
  \end{enumerate}
\end{lemma}

\begin{remark}[Failure of completeness for $\v_1$]
  \label{nocomp}
  \label{examrem}
  Lemma~\ref{exlem} shows that we can not have a completeness result
  (a converse of lemma~\ref{adeq} for closed terms) for $\v_1$.  To
  type the term $\l y^0.y^0$ by the type $(a \sqcap b) \f a$, we need
  an elimination rule for $\sqcap$ which we have in $\v_2$.  However,
  we will see that we have completeness for $\v_2$ if only one
  expansion variable is used.
\end{remark}

\section{Completeness of $\v_2$ with one expansion variable}
\label{complesec}

Let $a \in {\cal A}$, $e_1, e_2 \in {\cal E}$, $e_1 \not = e_2$ and
$Nat_0 = (e_1 a \f a) \f (e_2 a \f a)$.  Then:\\1) $\l f^0. f^0 \in
[Nat_0]$ and 2) It is not possible that $\l f^0. f^0 : \<() \v_2 Nat_0
\>$.

Hence $\l f^0. f^0 \in [ Nat_0]$ but $\l f^0. f^0$ is not
typable by $Nat_0$ and we do not have completeness in the presence
of more than one expansion variable.  The problem comes from the fact
that for the realisability semantics that we considered, we identify
all expansion variables. In order to give a completeness theorem we
will in what follows restrict our system to only one expansion
variable.  In the rest of this section, we assume that the set ${\cal
  E}$ contains only one expansion variable $e_c$.

The need of one single expansion variable is clear in part 2) of
lemma~\ref{inj} which would fail if we use more than one expansion
variable.  For example, if $e_1 \not = e_2$ then $e_1(e_2 a)^- = e_1 a
\not = e_2 a$. This lemma is crucial for the rest of this section and
hence, a single expansion variable is also crucial.

\begin{lemma}
  \label{inj}
  Let $U,V \in {\mathbb U}$ and $\deg(U) = \deg(V) > 0$.  1)
  $e_cU^- = U$ and 2) If $U^- = V^-$, then $U = V$.
\end{lemma}

Next, we divide  $\{y^n \mid y \in
{\cal V}_2\}$ disjointly amongst types of order $n$.

\begin{definition}
  Let $U \in {\mathbb U}$. We define sets of variables ${\mathbb V}_U$
  by induction on $\deg(U)$.  If $\deg(U) = 0$, then: ${\mathbb V}_U$
  is an infinite set of variables of degree $0$; if $y^0 \in {\mathbb
    V}_U$, then $y \in {\cal V}_2$; and if $U \neq V$ and $\deg(U) =
  \deg(V) = 0$, then ${\mathbb V}_U \cap {\mathbb V}_V = \emptyset$.
  If $\deg(U) = n+1$, then we put ${\mathbb V}_U = \{ y^{n+1} \mid y^n
  \in {\mathbb V}_{U^-}\}$.
\end{definition}

Our  partition of  ${\cal V}_2$ allows
useful infinite sets which contain type environments that
will play a crucial role in one particular type interpretation. These
sets and environments are given in the next definition.

\begin{definition}
  \label{infsetsenv}
  \begin{enumerate}
  \item Let $n \in {\mathbb N}$.  We let ${\mathbb G}^n=\{(y^n:U) \mid
    U\in{\mathbb U}$, $\deg(U) = n$ and $y^n \in {\mathbb V}_U \}$
    and ${\mathbb H}^n=\bigcup_{m\geq n}{\mathbb G}^m$. Note that
    ${\mathbb G}^n$ and ${\mathbb H}^n$ are not type environments
    because they are infinite sets.
  \item Let $n \in {\mathbb N}$, $M\in{\cal M}$ and $U\in{\mathbb U}$,
    we write $M : \< {\mathbb H}^n \v_2 U\>$ iff there is a type
    environment $\G \subset {\mathbb H}^n$ where $M : \< \G \v_2 U\>$
  \end{enumerate}
\end{definition}

Now, for every $n$, we define the set of the good terms of order $n$
which contain some free variable $x^i$ where $x \in {\cal V}_1$ and $i
\geq n$.

\begin{definition}
  Let $n \in {\mathbb N}$ and ${\cal V}^n = \{M \in {\mathbb
    M}^n \mid x^i \in FV(M)$ where $x \in {\cal V}_1$ and $i \geq n\}$.
  Obviously, if $n \in {\mathbb N}$ and
  $x \in {\cal V}_1$, then ${\cal N}^n_x \subseteq {\cal V}^n$.
\end{definition}

Here is the crucial  interpretation ${\mathbb I}$ for the proof of
completeness:

\begin{definition}
  Let ${\mathbb I}$ be the interpretation defined by:\\ for all type
  variables $a$,
  ${\mathbb I}(a) = {\cal V}^0 \cup \{M \in {\cal M}^0 \mid M : \<
  {\mathbb H}^0 \v_2 a\> \}$.
\end{definition}

${\mathbb I}$ is indeed an interpretation
and the interpretation of a type of order $n$ contains the
good terms of order $n$ which are typable in the special
environments which are parts of the infinite sets of
definition~\ref{infsetsenv}:

\begin{lemma}
  \label{comp}
  \begin{enumerate}
  \item ${\mathbb I}$ is an interpretation.  I.e., $\forall a \in
    {\cal A}$, ${\mathbb I}(a)$ is saturated and $\forall x \in {\cal
      V}_1$, ${\cal N}_x^0 \subseteq {\mathbb I}(a)\subseteq {\mathbb
      M}^0$.
  \item If $U\in{\mathbb U}$ is good and $\deg(U) = n$, then
    ${\mathbb I}(U) = {\cal V}^n \cup \{M \in {\mathbb M}^n \mid M : \<
    {\mathbb H}^n \v_2 U\> \}$.
  \end{enumerate}
\end{lemma}

${\mathbb I}$ is used to prove completeness (the proof is on the authors  web pages).

\begin{theorem}[Completeness]
  \label{comple}
  Let $U\in {\mathbb U}$ be good such that $\deg(U) = n$.
  \begin{enumerate}
  \item \label{compleone} $[U] = \{M \in {\mathbb M}^n \mid M : \< ()
    \v_2 U\> \}$.
  \item \label{completwo} $[U]$ is stable by reduction: i.e., if $M
    \in [U]$ and $M \rhd_\be^* N$, then $N \in [U]$.
  \item \label{complethree} $[U]$ is stable by expansion: i.e., if $N
    \in [U]$ and $M \rhd_\be^* N$, then $M \in [U]$.
  \end{enumerate}
\end{theorem}

\section{Conclusion and future work}
\label{concsection}

We studied the $\l I^{\mathbb N}$-calculus, an indexed version of the
$\l I$-calculus. This indexed version was typed using first an
intersection type system with expansion variables but without an
intersection elimination rule, and then using an intersection type
system with expansion variables and an elimination rule.

We gave a realisability semantics for both type systems showing that
the first type system is not complete in the sense that there are
types whose semantic meaning is not the set of $\l I^{\mathbb
  N}$-terms having this type.  In particular, we showed that $\l
y^0.y^0$ is in the semantic meaning of $(a \sqcap b) \f a$ but it is
not possible to give $\l y^0.y^0$ the type $(a \sqcap b) \f a$. The
main reason for the failure of completeness in the first system is
associated with the failure of the subject reduction property for this
first system.  We showed that the second system has the desirable
properties of subject reduction and expansion and strong normalisation
but that completeness fails if we use more than one expansion
variable. We then showed that completeness succeeds if we restrict the
system to one single expansion variable.

Because we show in the appendixes of the long version of this article
(which can be downloaded on the web page of the authors) that each of
these type systems, when restricted to the normal $\l I$-calculus
represents a well known intersection type system with expansion
variables, our study can be said to be the first denotational
semantics study of intersection type systems with expansion variables
(using realisability or any other approach) and outlines the
difficulties of doing so.  Although we have in this paper limited the
study to the $\l I$-calculus, future work will include extending this
work to the full $\l$-calculus and with an $\o$-type rule as well.

\bibliographystyle{jbwc}
{\bibliography{../included/literature,macros,bibliography,conferences}}

\end{document}